# DECISION MAKING: I I I - INCOMPLETE  INITIAL INFORMATION

## V. ZHUKOVIN,  Z. ALIMBARASHVILI

ABSTRACT: In the item there are presented Interval Structure allow to extend the class of multicriteria  decision making problems including  incomplete Initial Information. this step will facilitate dada gathering in pairs matching for decision making.

**Introduction. 1.** Some information about multicriteria   decision making problems. Given finite set $X$ of competitive decisions (alternatives). Its elements are   $x_i$, $i = 1,...,n$, that is   $x_i \in X$. A information structure is determinate on this set: Str. Then decision making problem may be presented by pair:

$$D = \langle X, Str \rangle \tag{1}$$

Information structure is formed on the bast of  data which we receive from different origin. It may be presented by different mathematical  means. Le1 examine two of them.

a) Given $m$ - dimensional criteria space of the base of   $m$ effectiveness criteria   $K_j$, $j = 1,...,m$ a point estimate (a point) correspond all   $x_i \in X$ in this space. Thus following vector - criterion define this information structure in this case:

$$Str^k : K(x_i) = \left\{ K_j(x_i) \right\}_{j=1}^m \text{ for all }  x_i \in X \tag{2}$$

Let us examine this structure. We have complete initial information. The   special preference relation $R_\Pi^k$ will be determined in pairs set      $E = X \bullet X$ which is named Pareto - domination. Pareto – set   $X_\Pi(R_\Pi^k) = X_\Pi^k$   correspond it. This is a set of effective decisions and only one of them may be chosen in decision making.

b) In this case information structure is presented as following:

$$Str^R : R = \left\{ R_j \right\}_{j=1}^m \tag{3}$$

When   $R_j \subset E$, $j = 1 \div m$  are scalar binary preference relations and $R$ is Vector Preference Relation (VPR).

Let us introduce the Pareto - domination for this case:

$$R_\Pi^R = \bigcap_{j=1}^m R_j \tag{4}$$

(4) and corresponding its the Pareto - set.

$$X_\Pi^R = X_\Pi(R_\Pi^R) \tag{5}$$

If all preference relation $R_j$, $j = 1 \div m$ are transitive then $R$ is transitive too and $X_{\Pi}^{R} \neq \varnothing$. It means that structure (2) is wider structure (2) because (2) is included in (2). This statement is verified by following transformation:

$$R_j = \left\{ (x_i, x_q) \in E \middle| K_j(x_i) \geq |K_j(x_q) \right\} \tag{6}$$

that is an effectiveness criterion $K_j$ always can be transformed in preference relation $R_j$ of course, Decision Maker (DM) must choose the acceptable compromise decision only in Paret – Set $X_{\Pi}^{R}$. This structure work only on the base of complete information, that is all decision pairs are ordered.

**Interval structure [2].** We isolate this class of decision making because it is principle in the description of the decision making problems with incomplete initial information. We notice it by IntStr. The intervals are closed and finite. They are situated in real-valued axis, that is:

$$\Delta = [a; b], \text{ where } b \geq a \tag{7}$$

(7) and $b$, $a$ are number from Re. Zero – interval corresponds $b = a$.

Let us examine some finite set of intervals:

$$In = \left\{ \Delta_j \right\}_{j=1}^{m} \tag{8}$$

(8) where $\Delta_j = [a_j; b_j]$ and $a_j, b_j \in \text{Re}$ for all $j = 1 \div m$. Let us introduce the binary preference relation $L^s$ defined in the pairs set $Q = In \times In$.

Interval Structure InStr is presented by following rule:

$$InStr : L^s = \left\{ (\Delta_i, \Delta_q) \middle| a_i \geq b_q \right\} \tag{9}$$

(9) Let us examine it. $L^s$ is strict(index $S$), unconnected (some pairs are mutually disjoint) and transitive. Therefore its Pareto-set is nonempty, that is we examine the interval set with finite $b = \max_{j=1 \div m} b_j$ and this is guarantee of exist of Pareto-set.

**I I I - Incomp1ecte Initial Information (variant $Str^k$).** Working with real objects very often we can't receive the point estimations of alternatives on the bases of vector = criterion. Then interval estimates may be used for description of incompleteness. In addition these intervals must following conditions:

- one interval correspond one point estimation of one $x_i \in X$ on one effectiveness criterion $K_j$ $j = 1 \div n$.

- in addition the point estimation always must be in corresponding it the interval.

- all intervals are finite, that is all intervals must be limited above below.

Thus when we receive and form incomplete initial information for decision making in fact we transform initial multicriteria decision making problem into the interval structure:

$$L = \left\{ L_j^s \right\}_{j=1}^{m} \tag{I0}$$

(10) Now this is with vector interval structure corresponding vector criterion. The following steps are known. Pareto - domination is defined:

$$L_{\Pi} = \bigcap_{j=1}^{m} L_j^s \tag{11}$$

(1l) This is a. sca1ar, strict unconnected and transitive preference relation defined on intervals set $\Lambda_0$, corresponding $X$ - set. On the base of this preference relation Pareto-set $X_{\Pi}^L$ formed. It is proved that $X_{\Pi}^k \subset X_{\Pi}^L$.

**I I I - Incomplete: Initial Information (variant $Str^R$).** principal element. It contain $m$ components, which are scalar binary preference relations $R_j$ $j = 1 \div m$. Let us examine one of them. Let us suppose that it is connected that is all pairs of decision are matching. In item [1] we presented the base of in this variant vector preference relation, noted VPR (2), is theory about Superiority Degree (SD) one decision over other. This is numerical skew - symmetry function $Z_j(x_i; x_q)$ defined on the pairs decisions that is on E-set. Having $R_j$ we can form corresponding SD. If $R_j$ is connected then we always can the base of SD form Utility Function (UF) defined on X- set and coordinated with $R_j$. Let note it by $\varphi_j(x_i)$ - order scale $Z(x_i; x_q)$ and $\varphi(x_i)$ may be determinate for VPR R(2) too.

Definition. VPR R[formula (2)] is unconnected even if one of $R_j$ $j = 1 \div m$ is unconnected .

If this case the incompleteness of initial information is connected with unconnected of VPR R. Let us give for $x_i \in X$ , $i = 1 \div n$, all decisions in $X$ which are incomparable with it.

$$x_i : \left\{ N_j(x_i) \right\}_{j=1}^m \qquad (12)$$

Where $N_j(x_i)$ is set of decisions in which are incomparable with $x_i$ epsilon X by $R_j$.

At first step we suppose that decision preference all $x_q \in N$ , where

$$N = N_1 \oplus N_2 \oplus ... \oplus N_j \oplus ... \oplus N_m \qquad (13)$$

The operant $\oplus$ consider that the intersection of these sets (let us say $r$ sets) is contained in $N$ - set r ones. This new information is added to initial information and now $R$ is connected. Then we can define upper utility function $\overline{\varphi}(x_i)$. Than we suppose that all $x_q \in N$ preference $x_i$ add this information to initial information. Now $R$ is connected too and we can define lower utility function $\underline{\varphi}(x_i)$, It holds:

$$\overline{\varphi}(x_i) \geq \underline{\varphi}(x_i) \quad \text{for all} \quad x_i \in X \qquad (14)$$

Where true utility function $\varphi(x_i)$ is in intervals:

$$\Delta_i = [a_i; b_i] , \ a_i = \underline{\varphi}(x_i) , \ b_i = \overline{\varphi}(x_i) \qquad (15)$$

Thus we transform the structure $Str^R$ in the structure. It is known what must be done later (the part 3 of this item)

**Additional Information (AI).** These structures must be realize on computer as the dialogue procedure of decision making with incomplete initial information (I I I) -procedure. Working with it in concrete problem situation we may receive additional information. In this item we don't examine false information but only true information. In contract initial intervals and form new interval structure noted by IntStr(a). Following results holds:

Statement.

1) For variant $Str^k$ :

$$X_{\Pi}^k \subseteq X_{\Pi}^L(a) \subseteq X_{\Pi}^L \qquad (16)$$

2) For variant $Str^R$ :

$$X_\Pi^R \subseteq X_\Pi^L(a) \subseteq X_\Pi^L \qquad (17)$$

where $X_\Pi^L(a)$ is Pareto -set of new interval structure, and $X_\Pi^k$ , $X_\Pi^R$ correspond complete information.

Thus we have a system (mathematical) of nonempty Pareto-sets, imbedded one to other. The similar result we have abstained for l -level preference relations [1;2]. On the base of it dialogue procedures effectively are formed because corresponding Pareto-set is connected with existing information.

**Inference. Conclusion.**

This structures must be realize on computer. In practice working with them in concrete situation we haven't possibility to obtain complete information at least initially. Our approach must be facilitate dale gathering - initial and additional. In this connection the interval structures play a important part.